\documentclass{amsart}
\usepackage{amssymb, latexsym}
\usepackage{amscd}

\newtheorem{theorem}{Theorem}

\newtheorem{remark}{Remark}
\newtheorem{definition}{Definition}

\begin{document}

\title{A Combination Theorem for $PD(n)$-Pairs}
\author{Rita Gitik}
\email{ritagtk@umich.edu}
\address{Department of Mathematics \\ University of Michigan \\ Ann Arbor, MI, 48109}  

\date{\today}

\begin{abstract} 
We prove a combination theorem for $PD(n)$-pairs.
\end{abstract}

\maketitle


\section{Introduction}

A group $G$ is an orientable $PD(n)$-group if for any coefficient $ZG$-module $M$ it has Poincar\'{e} duality, i.e. there exists a fundamental class $[G] \in H_n(G, \mathbf{Z})$ such that the cap product with $[G]$ induces an isomorphism 
$H^i(G, M) \overset{ [G]\cap}{\rightarrow} H_{n-i}(G, M)$, cf. \cite{D-D} p.136. 

A group $G$  and a family $\mathbf{S} = \{S_i, i \in I \}$ of not necessary distinct subgroups of $G$, denoted $(G, \mathbf{S})$, is called a group pair. The definition of an orientable $PD(n)$-group has been generalized to group pairs as follows.

\begin{definition} cf. \cite{B-E}, Definition 4.1.

A group pair $(G, \mathbf{S})$ is called a duality pair of dimension $n$ if there exists a $G$-module $C$, called the dualizing module of $(G, \mathbf{S})$, and an element $e \in H_n(G, \mathbf{S}; C)$
such that the cap-products $$H^k(G;A)\overset{e \cap} \rightarrow H_{n-k}(G, \mathbf{S}; C \otimes A)$$ and 
$$H^k(G, \mathbf{S}; A) \overset {e \cap} \rightarrow H_{n-k}(G;C \otimes A)$$ are isomorphisms for all $G$-modules $A$ and all $k \in \mathbf{Z}$.
\end{definition}

The dualizing module $C$ is isomorphic as a $G$-module to $H^n(G, \mathbf{S}; \mathbf{Z}G)) \cong H_0(G; C \otimes \mathbf{Z}G)$, cf. \cite{B-E}, p.299.

The paper \cite{Wa} has an excellent description of $PD(3)$-pairs.

\begin{definition} cf. \cite{B-E}, p. 302.
  
A duality pair $(G, \mathbf{S})$ of dimension $n$ is called a Poincar\'{e} duality pair of dimension $n$, in short a $PD(n)$-pair, if its dualizing module $C$ is isomorphic to $\mathbf{Z}$ as an Abelian group. In this case $C$ will be denoted $\tilde{\mathbf{Z}}$ and for any $G$-module $A$ the module $\tilde{\mathbf{Z}} \otimes A$ (with diagonal action) will be denoted $\tilde{A}$. If the $G$-action on $\tilde{\mathbf{Z}}$ is trivial, the $PD(n)$-pair is orientable, otherwise it is non-orientable.
\end{definition}

The main result of this paper is motivated by the topology of manifolds. 

Let $M_1$ and $M_2$ be compact $n$-manifolds with non-empty boundaries. Let $ \{ N^i_1, k \ge i \ge 1 \}$ and $\{ N^j_2, l \ge j \ge 1 \}$ be the boundary components of $M_1$ and $M_2$ respectively. Let $L_1$ and $L_2$ be compact submanifolds of $N^1_1$ and $N^1_2$ such that there exists a homeomorphism $f :(L_1, \partial L_1) \rightarrow  (L_2, \partial L_2)$. Then the space $M_1 \underset{f(L_1 )= L_2} \cup M_2$ is an $n$-manifold with the  boundary components $(N^1_1 - L_1) \underset{f(\partial L_1 )= \partial L_2} \cup (N^1_2 - L_2)$, $\{ N^i_1, k \ge i \ge 2 \}$, and  $\{ N^j_2, l \ge j > 2 \}$.

In order to formulate the corresponding statement for groups we need the following construction.

\begin{definition}
To simplify the exposition we consider a graph to be an oriented one-dimensional $CW$-complex.  Let $G_1$ and $G_2$ be the fundamental groups of graphs of groups $\Gamma_1$ and $\Gamma_2$ respectively. For any edge in those graphs which is a loop, insert an additional trivial vertex group  in its middle. Choose vertices $v_i$ in $\Gamma_i, i=1,2$ which have the same degree $n$ and let $\{ e_i^j, i=1,2, n \ge j \ge 1 \}$ be all the edges in $\Gamma_i$ beginning at $v_i$. Let $U_1$ and $U_2$ be the vertex groups of $v_1$ and $v_2$ respectively, and let $\{ P^i_j, i=1,2, n \ge j \ge 1 \}$  be the edge groups of  $\{ e_i^j, i=1,2, n \ge j \ge 1 \}$.  
Assume that there exists an isomorphism $\phi : (U_1, \overset {n} {\underset {j=1} \cup } P^j_1) \rightarrow (U_2, \overset {n} {\underset {j=1} \cup } P^j_2)$.

Define a new graph of groups, called the amalgamation of $\Gamma_1$ and $\Gamma_2$ along $\phi$, denoted by $\Gamma_1 \underset {\phi} * \Gamma_2$, as follows.

\begin{enumerate}
\item Form the disjoint union of $\Gamma_1$ and $\Gamma_2$.
\item Remove the vertices $v_i$ from $\Gamma_i$.
\item Identify each pair of the edges $ e_1^j$ and $ e_2^j, n \ge j \ge 1 $ by an orientation-reversing homeomorphism.
\end{enumerate}

Note that $\Gamma_1 \underset {\phi} * \Gamma_2$ need not be connected.
\end{definition}

The paper \cite{S-W} has an excellent discussion of graphs of groups.

\begin{theorem}
Let $(G_1, \mathbf{S_1} \cup (U_1 \underset {\partial U_1} * W_1))$ and $(G_2, \mathbf{S_2} \cup (U_2 \underset {\partial U_2} * W_2))$ be $PD(n)$-pairs such that $U_1 \underset {\partial U_1} * W_1$ and $U_2 \underset {\partial U_2} * W_2$ are fundamental groups of graphs of groups $\Gamma_1$ and $\Gamma_2$ respectively. Let $v_i, i=1,2$ be vertices in $\Gamma_i, i=1,2$ which have the same degree $n$ and let $\{ e_i^j, i=1,2, n \ge j \ge 1 \}$ be all the edges in $\Gamma_i$ beginning at $v_i$. Let $U_1$ and $U_2$ be the vertex groups of $v_1$ and $v_2$ respectively, and let $\{ P^i_j, i=1,2, n \ge j \ge 1 \}$  be the edge groups of  $\{ e_i^j, i=1,2, n \ge j \ge 1 \}$.  
Assume that there exists an isomorphism $\phi : (U_1, \overset {n} {\underset {j=1} \cup } P^j_1) \rightarrow (U_2, \overset {n} {\underset {j=1} \cup } P^j_2)$. Let $G$ be the amalgamated free product $G_1 \underset{\phi(U_1) = U_2}{\ast} G_2$. Note that the group $W_1 \underset {\phi(\partial U_1)= \partial U_2} * W_2$ is the set of the fundamental groups of the components of $\Gamma_1 \underset {\phi} * \Gamma_2$. If all the subgroups $\{ P^i_j, i=1,2, n \ge j \ge 1 \}$ are $PD(n-2)$-groups then $(G, ( W_1 \underset {\phi(\partial U_1)= \partial U_2} * W_2 )\cup \mathbf{S_1} \cup \mathbf{S_2})$ is a $PD(n)$-pair.
\end{theorem}

\begin{remark}
Bieri and Eckmann proved a special case of Theorem 1 for $(U_i \underset {\partial U_i} * W_i)=U_i, i=1,2$, cf. \cite{B-E}, Theorem 8.1. The main technical tool of their proof is Theorem 3.5 of \cite{B-E} which has the following topological interpretation. Let $X$ and $Y$ be compact $n$-manifolds such that the boundary of $X$ is $\partial X = \partial_1 X \cup T_X$ and the boundary of $Y$ is $\partial Y = \partial_1 Y \cup T_Y$ where $T_X$ and $T_Y$ are homeomorphic compact $(n-1)$-manifolds. Let $M= X \underset {T_X = T_Y} \cup Y$ and let $T$ be the image of $T_X$ (which is also the image of $T_Y$) in $M$. Note that the inclusions $$(M,\partial M) \overset i {\rightarrow} (M,\partial M \cup Y) \overset j {\leftarrow} (X,\partial X)$$ induce homomorphisms on the homology $H_k(M,\partial M) \overset {i_*} {\rightarrow} H_k(M,\partial M \cup Y) \overset {j_*} {\leftarrow} H_k(X,\partial X)$, where the map $j_*$ is the excision isomorphism. Hence there exists a map $\phi_X : H_k(M,\partial M) \rightarrow H_k (X,\partial X)$. Similarly  there exists a map $\phi_Y : H_k(M,\partial M) \rightarrow H_k (Y,\partial Y)$. The topological analog of Theorem 3.5 of \cite{B-E} is that the following sequence $$ \rightarrow H_k(T) \rightarrow H_k(M,\partial M) \overset {(\phi_X,- \phi_Y)} {\rightarrow} H_k(X,\partial X)\oplus H_k(Y, \partial Y) \rightarrow H_{k-1}(T) \rightarrow$$ is exact. When $k=n$ the first group in that sequence is trivial,  the second and the fourth groups in that sequence are isomorphic to $\mathbf{Z}$, and the third group is isomorphic to $\mathbf{Z \oplus Z}$.
\end{remark}

In order to prove Theorem 1 we will need the following technical result. 

\begin{theorem}
Let $(G_i, \mathbf{S_i}), i=1,2$ be group pairs of dimension $n$. Let $U_i \underset {\partial U_i} * W_i$ be a subgroup of $G_i, i=1,2$ such that 
$(U_1, \partial U_1)$ is isomorphic to $(U_2,\partial U_i)$. Let $\overline{\mathbf{S_i}}= \mathbf{S_i} \cup W_i, i=1,2$ and let $G=G_1 \underset {U_1=U_2}
* G_2$. Denote the image of $(U_1, \partial U_1)$ in $G$ (which is equal to the image of $(U_2, \partial U_2)$ in $G$) by $(U, \partial U)$.
Denote $\overline{\mathbf{S}} = \overline{\mathbf{S_1}} \cup \overline{\mathbf{S_2}} \cup W_1 \cup W_2$. Denote $\partial G = \mathbf{S_1} \cup \mathbf{S_2} \cup (W_1 \underset {\partial U_1 = \partial U_2} * W_2)$ and for $i=1,2$ denote $\partial G_i = \mathbf{S_i}  \cup (U_i \underset {\partial U_i} * W_i)$. Then the following long sequence $$ \rightarrow H_k(U, \partial U) \rightarrow H_k(G, \partial G) \rightarrow H_k(G_1, \partial G_1) \oplus H_k(G_2, \partial G_2) \rightarrow H_{k-1}(U, \partial U) \rightarrow $$ is exact. 
\end{theorem}

\section{Proof of Theorem 2}

Consider the following commutative diagram
\[
\begin{CD}
H_k(\partial U)   @>0>>  H_{k-1}(\partial U) @>(+1,-1)>> H_{k-1}(\partial U) \oplus H_{k-1}(\partial U) @>>> H_{k-1}(\partial U)\\
@VVV                 @VVV                                         @VVV                           @VVV\\
H_k(U) @>>> H_k(G, \overline{\mathbf{S}}) @>>> H_k(G_1, \overline{\mathbf{S_1}} \cup U) \oplus H_k(G_2,  \overline{\mathbf{S_2}} \cup U) @>>> H_{k-1}(U)\\
@VVV                 @VVV                                         @VVV                           @VVV\\
H_k(U, \partial U) @>>> H_k(G, \partial G) @>>> H_k(G_1, \partial G_1) \oplus H_k(G_2, \partial G_2) @>>> H_{k-1}(U, \partial U)\\
@VVV                 @VVV                                         @VVV                           @VVV\\
H_{k-1}(\partial U)   @>0>>  H_{k-2}(\partial U) @>(+1,-1)>> H_{k-2}(\partial U) \oplus H_{k-2}(\partial U) @>>> H_{k-2}(\partial U)\\
\end{CD}
\]

Theorem 3.5 of \cite{B-E} implies that the second row of this diagram is exact. The first and the last columns of this diagram are exact because they are long exact sequences of pairs. The first and the last rows of the diagram split into short exact sequences.

Compare $H_k(G, \overline{\mathbf{S}})$ and $H_k(G, \partial G)$ using the exact sequence of the triple $(G, \partial G, \overline{\mathbf{S}})$. $$H_k(\partial G, \overline{\mathbf{S}}) \rightarrow H_k(G, \overline{\mathbf{S}}) \rightarrow H_k(G, \partial G) \rightarrow H_{k-1}(\partial G, \overline{\mathbf{S}})$$ where the first term is, by definition, $$H_k(\mathbf{S_1} \cup \mathbf{S_2} \cup (W_1 \underset {\partial U} * W_2), \mathbf{S_1} \cup \mathbf{S_2} \cup W_1 \cup W_2)$$ which is isomorphic by excision to $H_k(\partial U \times I, \partial U \times \partial I) \cong H_{k-1}(\partial U)$. It follows that the second column of the diagram is exact.
 
Compare $H_k(G_1, \overline{\mathbf{S_1}} \cup U)$ and $H_k(G_1, \partial G_1)$ using the exact sequence of the triple $(G_1, \partial G_1, \overline{\mathbf{S_1}} \cup U)$. $$H_k(\partial G_1, \overline{\mathbf{S_1}} \cup U) \rightarrow H_k(G_1, \overline{\mathbf{S_1}} \cup U) \rightarrow H_k(G_1, \partial G_1) \rightarrow H_{k-1}(\partial G_1, \overline{\mathbf{S_1}} \cup U)$$ where the first term  is, by definition,  $$H_k(\mathbf{S_1} \cup ( W_1 \underset {\partial U} * U), \mathbf{S_1} \cup W_1 \cup U)$$
 which is isomorphic by excision to $H_k(\partial U \times I, \partial U \times \partial I) \cong H_{k-1}(\partial U)$. It follows that the third column of the diagram is exact. 
 
Using diagram chasing, it is easy to show that the third row of the diagram is exact, completing the proof of Theorem 2. 
 
 Note that if $G$ is a $PD(k)$-group then $(G, \{G, G \})$ is a $PD(k+1)$-pair.

\section{Proof of Theorem 1}

We use the notation of Theorem 2.

For $i=1,2$ let $\Gamma'_i$ be the the graph $\Gamma_i - v_i- (\overset {n} {\underset{j=1} \cup } e_i^j)$, and let $W_i$ be the fundamental group of $\Gamma'_i$. Note that $\Gamma'_i$ need not be connected. Hence for $i=1,2$ the group $U_i \underset {\partial U_i} * W_i$ can be considered as the fundamental group of the graph of groups with two vertex groups, namely $U_i$ and $W_i$, and $n$ edge groups, namely $\overset {n} {\underset{j=1} \cup } P^j_i$. By assumption,
$U_i \underset {\partial U_i} * W_i$ is a $PD(n-1)$-group. Also by assumption, the subgroups $\{ P^i_j, i=1,2, n \ge j \ge 1 \}$ are $PD(n-2)$-groups, hence part (ii) of Theorem 8.4 of \cite{B-E} implies that $(U_1,\overset {n} {\underset{j=1} \cup } P^j_1)$ and $(U_2, \overset {n} {\underset{j=1} \cup } P^j_2)$ are $PD(n-1)$-pairs.

For $i=1,2$ Theorem 8.4 of \cite{B-E} also implies that the dualizing modules of $U_i \underset {\partial U_i} * W_i$ and $(U_i,\overset {n} {\underset{j=1} \cup } P^j_i)$ are isomorphic to $H^n((U_i \underset {\partial U_i} * W_i), \mathbf{Z}(U_i \underset {\partial U_i} * W_i))$, considered as $U_i$-modules by restriction. Denote this dualizing module $C$. It has a $G$-module structure.

Theorem 2 states that the sequence 
$$0 \rightarrow H_n(U, \partial U) \rightarrow H_n(G, \partial G) \rightarrow H_n(G_1, \partial G_1) \oplus H_n(G_2, \partial G_2) \rightarrow H_{n-1}(U, \partial U)\rightarrow $$
is exact.

Choose fundamental classes $e_i \in  H_n(G_i,  \partial G_1), i=1,2$ and $e_0 \in H_{n-1}(U, \partial U)$ such that $\partial e_1 = \partial e_2 = e_0$, and consider the element $e \in H_n(G, \partial G)$ which maps onto $(e_1, -e_2)$. The diagram in the proof of Theorem 2 shows that $(G, \partial G)$ is a $PD(n)$-pair with dualizing module $C$ and fundamental class $e$, completing the proof of Theorem 1.

\section{Acknowledgment}

We would like to thank Peter Scott for helpful discussions.

\end{document}